\documentclass{amsart}

\RequirePackage[left,mathlines]{lineno}%
\RequirePackage{amsthm}%
\RequirePackage{amsmath}%
\RequirePackage{amsxtra}%
\RequirePackage{amstext}%
\RequirePackage{amssymb}%
\RequirePackage{latexsym}%
\RequirePackage{dsfont}%

\theoremstyle{definition}

\theoremstyle{remark}

\numberwithin{equation}{section}

\ifx\pdfoutput\undefined%
   \IfFileExists{graphicx.sty}{\RequirePackage[dvips]{graphicx}
   \DeclareGraphicsExtensions{.eps,.ps}}{}
\else%
   \ifnum\pdfoutput=0
      \IfFileExists{graphicx.sty}{\RequirePackage[dvips]{graphicx}
      \DeclareGraphicsExtensions{.eps,.ps}}{}
   \else
      \IfFileExists{graphicx.sty}{\RequirePackage[pdftex]{graphicx}
      \DeclareGraphicsExtensions{.pdf,.png,.jpg}}{}
   \fi
\fi 

\RequirePackage[numbers]{natbib}

\RequirePackage[colorlinks,citecolor=blue,urlcolor=blue,pagebackref]{hyperref}
\RequirePackage{hypernat}%
\RequirePackage{color}%

\definecolor{red}       {rgb}{0.0,0.0,1.0}        
\definecolor{magenta}   {rgb}{0.0,0.0,1.0}        
\definecolor{cyan}      {rgb}{0.0,0.0,1.0}        
\definecolor{green}     {rgb}{0.0,0.4,0.3}        

\begin{document}



\title[Funny Problems in Intuitive Topology]%
{Funny Problems in Intuitive Topology}


\author{R. Tavakoli}

\address{Department of Material Science and Engineering, Sharif University of
Technology, Tehran, Iran, P.O. Box 11365-9466}

\thanks{\\Department of Material Science and Engineering, Sharif University of
Technology, Tehran, Iran, P.O. Box 11365-9466,
\href{mailto:tav@mehr.sharif.edu}{tav@mehr.sharif.edu},
\href{mailto:rohtav@gmail.com}{rohtav@gmail.com}.}

\email{\href{mailto:tav@mehr.sharif.edu}{tav@mehr.sharif.edu},
\href{mailto:rohtav@gmail.com}{rohtav@gmail.com}}




\maketitle%



\begin{abstract}%

The goal of this article is to introduce some beautiful known
riddles in intuitive topology; hoping to make at least some fun for
the reader.

\end{abstract}%


\section{Funny problems in intuitive topology}

The materials of this section are adapted from reference:
\cite{prasolov-intuitive}.

\begin{figure}[ht]%
\begin{center}%
\includegraphics[width=8.cm]{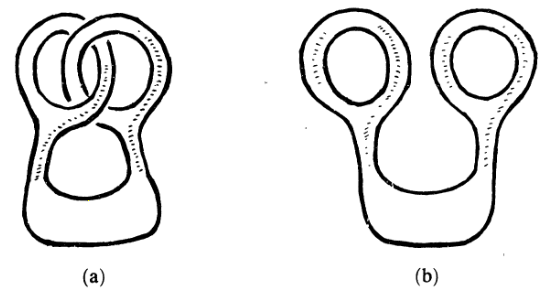}%
\caption{Plot related to statement of problem 1 \cite{prasolov-intuitive}.}%
\label{fig:p1}%
\end{center}%
\end{figure}

Problem 1. (\cite{prasolov-intuitive}) Show that the elastic body
represented in figure \ref{fig:p1} (a) can be deformed so as to
become the one shown in \ref{fig:p1} (b). In other words, were the
human body elastic enough, after making linked rings with your index
fingers and thumbs, you could move your hands apart without
separating the joined fingertips.

\begin{figure}[ht]%
\begin{center}%
\includegraphics[width=8.cm]{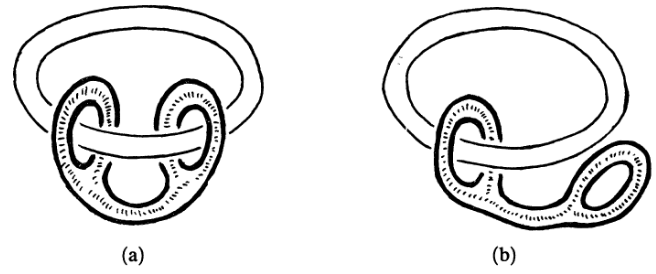}%
\caption{Plot related to statement of problem 2 \cite{prasolov-intuitive}.}%
\label{fig:p2}%
\end{center}%
\end{figure}

Problem 2. (\cite{prasolov-intuitive}) A pretzel has two holes that
"hold" a doughnut (see figure \ref{fig:p2} (a)). Show that the
pretzel can be deformed in such a way that one of its "handles" will
unlink itself from the doughnut (figure \ref{fig:p2} (b)).

\begin{figure}[ht]%
\begin{center}%
\includegraphics[width=8.cm]{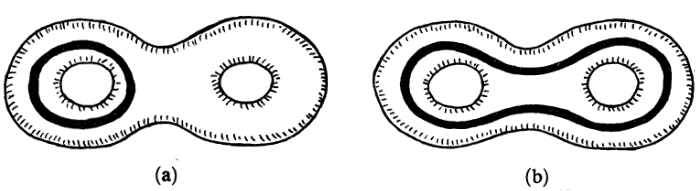}%
\caption{Plot related to statement of problem 3 \cite{prasolov-intuitive}.}%
\label{fig:p3}%
\end{center}%
\end{figure}

Problem 3. (\cite{prasolov-intuitive}) A circle is drawn on a
pretzel with two holes (\ref{fig:p3} (a)). Show that it is possible
to deform the pretzel so that the circle will be in the position
represented in \ref{fig:p3} (b).

Problem 4. (\cite{prasolov-intuitive}) Show that a punctured tube
from a bicycle tire can be turned inside out. More precisely, this
would be possible if the rubber from which the tube is made were
elastic enough. In real life it is impossible to turn a punctured
tube inside out.

\begin{figure}[ht]%
\begin{center}%
\includegraphics[width=8.cm]{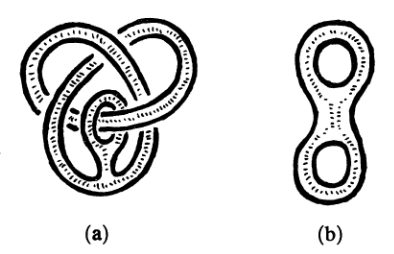}%
\caption{Plot related to statement of problem 5 \cite{prasolov-intuitive}.}%
\label{fig:p5}%
\end{center}%
\end{figure}

Problem 5. (\cite{prasolov-intuitive}) Show that the fancy pretzel
represented in \ref{fig:p5} (a) can be deformed into the ordinary
pretzel with two holes (\ref{fig:p5} (b)).

\noindent \textcolor[rgb]{0.00,0.00,1.00}{{\bf{ Do you can resolve
these problems by yourself (please try and then go to the next
page)?}}}

\clearpage

\begin{figure}[ht]%
\begin{center}%
\includegraphics[width=8.cm]{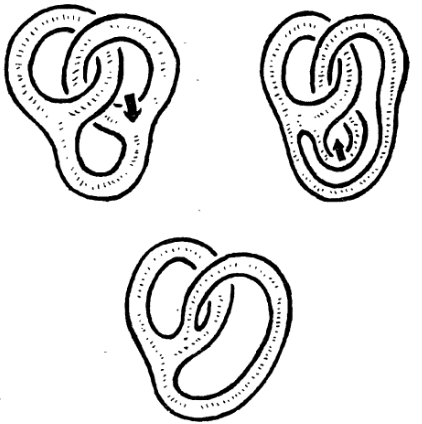}%
\caption{graphical solution for problem 1 \cite{prasolov-intuitive}.}%
\label{fig:sp1}%
\end{center}%
\end{figure}

\begin{figure}[ht]%
\begin{center}%
\includegraphics[width=8.cm]{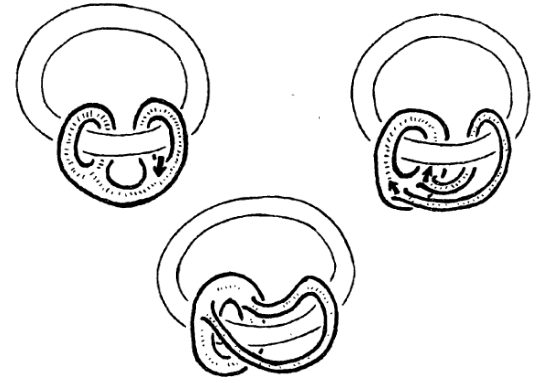}%
\caption{graphical solution for problem 2 \cite{prasolov-intuitive}.}%
\label{fig:sp2}%
\end{center}%
\end{figure}

\begin{figure}[ht]%
\begin{center}%
\includegraphics[width=12.cm]{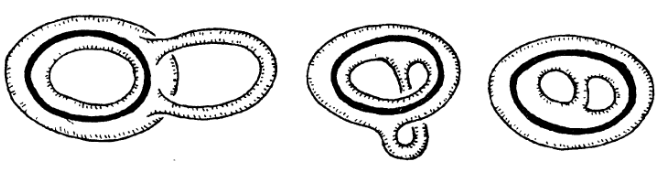}%
\caption{graphical solution for problem 3 \cite{prasolov-intuitive}.}%
\label{fig:sp3}%
\end{center}%
\end{figure}

\begin{figure}[ht]%
\begin{center}%
\includegraphics[width=12.cm]{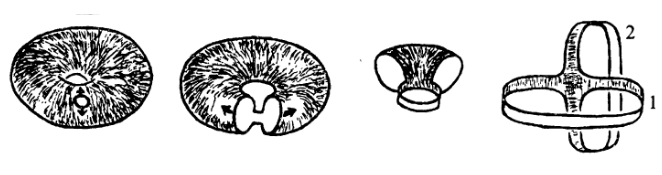}%
\caption{graphical solution for problem 4 \cite{prasolov-intuitive}.}%
\label{fig:sp4}%
\end{center}%
\end{figure}

\begin{figure}[ht]%
\begin{center}%
\includegraphics[width=8.cm]{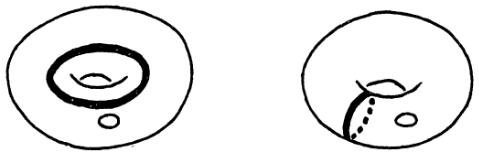}%
\caption{graphical solution for problem 4 (continue) \cite{prasolov-intuitive}.}%
\label{fig:sp4_2}%
\end{center}%
\end{figure}

\begin{figure}[ht]%
\begin{center}%
\includegraphics[width=14.cm]{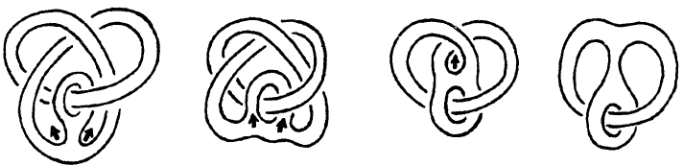}%
\caption{graphical solution for problem 5 \cite{prasolov-intuitive}.}%
\label{fig:sp5}%
\end{center}%
\end{figure}

Graphical solutions to above problems are introduced in figures
\ref{fig:sp1} through \ref{fig:sp5}.


Regarding to the problem 4, first we perform the deformations shown
in figure \ref{fig:sp4}. Then we can change the position of the
obtained figure so that its "inside" (shown in white) becomes its
 "outside" (the shaded side of the surface) and vice-versa, simply by
moving it as a rigid body in space until the hoop 1 occupies the
position of the hoop 2. Once this is done, the previous deformations
performed in reverse order result in the tube being turned inside
out as required. Note that this procedure interchanges the
"parallel" and the "meridian" of the tube (see figure
\ref{fig:sp4_2}) \cite{prasolov-intuitive}.

Regarding to the problem 5, first we perform the deformation shown
in figure \ref{fig:sp5}. The solid thus obtained (provided it is
elastic) can clearly be def.ermed into the one shown in figure
\ref{fig:p1} (a). It now remains to apply the solution of Problem 1
(cf. figure \ref{fig:sp1})

Topologists call such objects isotopic (in space). More clearly,
when we consider objects in space "up to deformations", i.e., we did
not distinguish objects that can be deformed into each other by
reshaping them, we called such objects as isotopic.

In the same way, the topology aid us to identify isotopic knots. For
example all knots represented on the top row of figure \ref{pk1} can
be deformed into each other. Similarly, all knots represented on the
bottom row of figure \ref{pk1} can be deformed into each other. For
solution readers are refereed to chapter 2 of
\cite{prasolov-intuitive}.

\begin{figure}[ht]%
\begin{center}%
\includegraphics[width=14.cm]{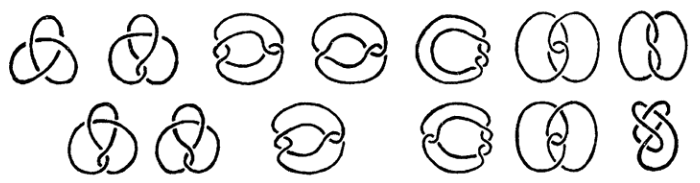}%
\caption{All knots represented on the top row are isotopic,
and the same for the bottom row \cite{prasolov-intuitive}.}%
\label{fig:pk1}%
\end{center}%
\end{figure}%
%


\bibliographystyle{unsrt} 
\bibliography{biblio}%

\end{document}